\documentclass{amsart}

\usepackage{amsmath,amsthm,amsopn,amstext,amscd,amsfonts,amssymb}
\usepackage{dsfont, fontenc} % para escribir correctamente los cjtos R,Z,C...
\usepackage{comment}
\usepackage[active]{srcltx}
\usepackage{graphicx}
\usepackage{caption}
\usepackage{subcaption}

\usepackage{algpseudocode}
\newtheorem*{problemq}{\sc Question (Q)}

\newtheorem{theorem}{\sc Theorem}[section]
\newtheorem{obs}{\sc Observation}

\usepackage[active]{srcltx}
\usepackage{algorithm}

\makeatletter
\def\BState{\State\hskip-\ALG@thistlm}
\makeatother

\def\downbar#1{
\setbox10=\hbox{$#1$}
            \dimen10=\ht10 \advance\dimen10 by 2.5pt
            \ifdim \dimen10<15pt %equals approximately 0.5cm
               \advance\dimen10 by -0.5pt
               \dimen11=\dimen10
               \advance\dimen10 by 2.5pt
               \lower \dimen11
            \else \lower \ht10 \fi
            \hbox {\hskip 1.5pt \vrule height \dimen10 depth \dp10}}
\def\upbar#1{
\setbox10=\hbox{$#1$}
            \dimen10=\ht10 \advance\dimen10 by \dp10 \advance\dimen10 by 2.5pt
            \ifdim \dimen10<15pt %equals approximately 0.5cm
                \advance\dimen10 by 2pt \fi
            \raise 2.5pt \hbox {\hskip -1.5pt \vrule height \dimen10}}

\newcommand{\re}{\mathbb{R}}

\newcommand{\dps}{\displaystyle}

\makeatletter\renewcommand\theenumi{\@roman\c@enumi}\makeatother

\begin{document}
\title[On zeros of polynomials in best $L^p$-approximation]{On zeros of polynomials in best $L^p$-approximation and inserting mass points}
\author{K. Castillo}
\address{CMUC, Department of Mathematics, University of Coimbra,  3001-501 Coimbra, Portugal}
\email{kenier@mat.uc.pt}
\author{M. S. Costa}
\author{F. R. Rafaeli}
\address{FAMAT-UFU, Department of Mathematics, University of Uberl\^andia, 38408-100 Uber{l\^a}ndia, Minas Gerais, Brazil}
\email{marisasc@ufu.br}
\email{rafaeli@ufu.br}

\subjclass[2010]{30C15}
\date{\today}
\keywords{Polynomials, minimal $L^p$ norm, monotonicity, zeros}
\begin{abstract}
The purpose of this note is to revive in $L^p$ spaces the original A. Markov ideas to study monotonicity of zeros of orthogonal polynomials. This allows us to prove and improve in a simple and unified way our previous result [Electron. Trans. Numer. Anal., 44 (2015), pp. 271--280] concerning the discrete version of A. Markov's theorem on monotonicity of zeros.\end{abstract}
\maketitle

\section{Introduction and main results}
Let $\mu$ be a positive and nontrivial Radon measure on a compact set $A \subset \re$. For $1 < p < \infty$, the space $L^p(\mu)$ denotes the set of all equivalent classes of $\mu$-measurable functions $f$ such that $|f|^p$ is $\mu$-summable, endowed with the usual vector operations and with the norm
\begin{align}\label{norm}
\| f \|_{p}:=\left(\int |f(x)|^{p}d\mu(x)\right)^{1/p}.  
\end{align}
Set $X:=L^p(\mu)$. By a well known result by Clarkson \cite[Corollary, p. $403$]{C36}, $X$ is uniformly convex. Following Bourbaki \cite[Definition {\rm I}, p. $166$]{B68}, define $\mathbb{N}:=\{0,1,\dots\}$. Fix $n \in \mathbb{N}$ and set $K:=\mathcal{P}_{n}$,  $\mathcal{P}_{n}$ being the set of all real polynomials of degree at most $n$ regarded as a subspace of $X$. Since $K$ is finite dimensional, $K$ is a closed convex subspace of $X$. %Following Singer \cite[p. $15$]{S70}, $\mathfrak{L}_K(f)$ denotes the set of all elements of best approximation of $f\in X$ by elements of $K$. 
It is known that for any point $f \in X$, there is a unique point $g_0 \in \mathfrak{L}_K(f)$ (cf. \cite[Theorem $8$, p. $45$]{L02}). The preceding affirmation thus guarantees the existence and uniqueness of $g_0 \in \mathfrak{L}_K(x^{n+1})$. By the characterization of elements of  best approximation (cf. \cite[Theorem $1.11$]{S70}) $g_0 \in \mathfrak{L}_K(x^{n+1})$ if and only if
\begin{align}\label{char}
\int g(x) |x^{n+1}-g_0(x)|^{p-1} \text{sgn}(x^{n+1}-g_0(x))d\mu(x)=0 \quad (g\in K).
\end{align}
Consider the (monic) polynomial $P_{n+1,p}(x):=x^{n+1}-g_0(x)$. As a consequence of \eqref{char}, the minimum of the norm \eqref{norm} taken over all (monic)  real polynomials $P_{n+1}$ of degree $n+1$ is attained when $P_{n+1}:=P_{n+1,p}$. By Fej\'er's convex hull theorem (cf. \cite[Theorem $10.2.2$]{D75}), the zeros of $P_{n+1,p}$ all lie in the closure of the convex hull of $\text{supp}(\mu)$. Furthermore, all the zeros of $P_{n+1,p}$ are simple \footnote{Suppose, contrary to our claim, that $x_0$ is a multiple zero. From \eqref{char} we have
$$
\int \frac{P_{n+1,p}(x)}{(x-x_0)^2}|P_{n+1,p}(x)|^{p-1} \text{sgn}(P_{n+1,p}(x))d\mu(x)=\int \frac{|P_{n+1,p}(x)|^{p}}{|x-x_0|^2}d\mu(x)=0,
$$
a contradiction.}.

The central concern of this work is the following

\begin{problemq}\label{P}%Set $C:=$. such that \eqref{moments} holds, {\em mutantis mutandis}
Let $\mu$ be a positive and nontrivial Radon measure on a compact set $A \subset \re$. Assume that $\mathrm{d}\mu(x,t)$ has the form \footnote{The Dirac measure $\delta_{y}$ is a positive Radon measure whose support is the set $\{y\}$.}
\begin{align}\label{measure}
\mathrm{d}\alpha(x,t)+\jmath(t) \delta_{y(t)},
\end{align}
where $\mathrm{d}\alpha(x,t):=\omega(x,t)\mathrm{d}\nu(x)$ ($\omega$ is a positive weight  and $\nu$ is a positive and nontrivial Radon measure) and, $\jmath(t) \in \re_+$ and $y(t)\in \re$ are continuous differentiable function of $t \in U$, $U$ being an open interval on $\re$. Determine sufficient conditions in order for the zeros of the polynomial $P_{n+1,p}(x,t)$ $(2 \leq p<\infty)$ to be strictly increasing functions of $t$.
\end{problemq}

For reasons of economy of exposition, we intentionally avoided the case $1<p<2$ or the one in which we have infinitely many mass points. Even though the reader has to proceed with caution in these cases, under natural additional assumptions, Theorem \ref{main} below remains true, {\em mutatis mutandis}.  In \cite{CCR2} the reader can find a detailed study of the case $p=2$ when we have infinitely many mass points, using the ideas originally presented in this work. We recall that this solves an open problem posed by Ismail at the end of the $1980$'s within the framework of orthogonal polynomials (cf. \cite[Problem $1$]{I89} and \cite[Problem $24.9.1$]{I05}). When \eqref{measure} has the form $\omega(x,t)\mathrm{d}x$ and $p=2$, Question (Q) was studied as early as $1886$ by A. A. Markov \cite[p. $178$]{M86}, in a work with many lights and some shadows (see, for instance, \cite[Section $1$]{C17} for some historical remarks). When \eqref{measure} has the form $\omega(x,t)\mathrm{d}\nu(x)$ and $p=2$, Question (Q) was posed as an exercise in Freud's book \cite[Problem $16$, p. $133$]{G71} (a proof of such result can be found in the more recent book by Ismail \cite[Theorem $7.1.1$]{I05}). When \eqref{measure} has the form $\omega(x,t)\mathrm{d}x$, $A:=[-1,1]$, and $1\leq p \leq \infty$,  Question (Q) was studied by Kro\'o and Peherstorfer \cite{KP87}. When \eqref{measure} has the form $\omega(x)\mathrm{d}x+\jmath \delta_{y(t)}$ and $p=2$, Question (Q) was considered in \cite[Theorem $2.2$]{CR15} through a combination of elementary facts. It is, therefore, natural that this last result be broadened to $L^p$ spaces. Not surprisingly, this can be easily achieved by using Markov's original ideas \footnote{In his classical book \cite[Footnote $31$, p. $116$]{S75}, Szeg\H{o} refers his proof of Markov's theorem in the following terms: {\em ``This proof does not differ essentially from the original one by Markov, although the present arrangement is somewhat clearer."}. Probably this assertion has avoided the attention of some mathematicians to Markov's work. While it is true that in the framework of orthogonal polynomials Szeg\H{o}'s argument becomes especially elegant, Markov's approach works in a more general framework. Szeg\H{o}'s approach is based on Gauss mechanical quadrature, which was an approach that Stieltjes suggested to handle the problem, see \cite[Section $5$, p. $391$]{S87}.}. Our main result reads as follows:
\begin{theorem}\label{main}
Assume the notation and conditions of Question (Q). Assume further the existence and continuity for each $x \in A$ and $t \in U$ of $(\partial \omega/\partial t) (x,t)$. % as well as the convergence of the integrals
%$$
%\int x^m \frac{\partial \omega}{\partial t}(x,t)d\nu(x) \quad (m=0,\dots, n \,p-1),
%$$
%uniformly in every closed subset of $U$. 
Denote by $x_0(t), \dots, x_n(t)$ the zeros of $P_{n+1,p}(x,t)$. Fix $k \in \{0,\dots, n\}$ and set
 $$
 d_k(t):=
  \begin{cases}
    y(t)-x_k(t)      & \quad \text{if } y(t)\not = x_k(t),\\
    1  & \quad \text{if }  y(t)= x_k(t).\\
  \end{cases}
 $$
 Define the function
 $$
 R_k(t):=\sum_{j=0}^n{\vphantom{\sum}}' \frac{p-\delta_{j,k}}{y(t)-x_j(t)},
 $$
where the prime means that the sum is over all values $j$ and $t$ for which $ y(t)\not = x_j(t)$. Then $(\mathrm{d}x_k/\mathrm{d}t)(t)$ is strictly positive for those values of $t$ such that
\begin{align}\label{eq}
\frac{1}{d_k(t)} \left\{ \frac{\jmath'(t)}{\jmath(t)}+ y'(t) R_k(t)-\dps \frac{1}{\omega(x_k(t),t)}\dps\frac{\partial \omega}{\partial t}(x_k(t),t)\right\}&\geq 0,
\end{align}
and
\begin{align}\label{eq2}
\frac{1}{\omega(x,t)}\frac{\partial \omega}{\partial t}(x,t)
\end{align}
is an increasing function of $x \in A$, provided that at least the inequality \eqref{eq} be strict or the function \eqref{eq2} be nonconstant on $A$.
\end{theorem}

 The next observations concern the cases studied in the literature for $p=2$. As far as we know, these are the only ones that have been studied up to now. It is worth highlighting that such cases are the simplest consequences that can be derived from Theorem \ref{main}.

\begin{obs}\label{coro1}\footnote{Observation \ref{coro1} for $p=2$ was proved for the first time in \cite[Theorem $2.2$]{CR15}. In order to have monotonicity of zero the location of the mass point outside $\rm{Co}(A)$ is quite natural. In this regard, the statements of Theorem $2$ and Corollary $3$ in {\tt arXiv:1501.07235 [math.CA]} appear to be incorrect.}
Assume the notation and conditions of Theorem \ref{main} under the constraint that  $\mathrm{d}\mu(x,t)$ $=\mathrm{d}\alpha(x)+\jmath \delta_{y(t)}$. Define the sets \footnote{$A^c:=\{x\in \re \ | \ x\not \in A\}$ and $\rm{Co}(A)$ denotes the convex hull of $A$.}
\begin{align*}
B_-&:=\{t \in U \ | \ y(t)\in \rm{Co}(A)^c \wedge y'(t)<0\},\\
B_+&:=\{t \in U \ | \ y(t)\in \rm{Co}(A)^c \wedge y'(t)>0\}.
\end{align*}
Then all the zeros of $P_{n+1,p}(x,t)$ are strictly decreasing (respectively, increasing) functions of $t$ on $B_-$ (respectively, on $B_+$).
\end{obs}

\begin{obs} \footnote{The case $p=2$, often considered in the literature, can be easily handled by using very elementary results.}
Assume the notation and conditions of Theorem \ref{main} under the constraint that $\mathrm{d}\mu(x,t)$ $=\mathrm{d}\alpha(x)+\jmath(t) \delta_{y}$. Define the sets
\begin{align*}
C_-:=\{t \in U \ | \  \jmath'(t)<0\}, \quad C_+:=\{t \in U \ | \ \jmath'(t)>0\}.
\end{align*}
If $x_k(t)<y$ (respectively, $x_k(t)>y$) for each $t\in U$, then $x_k(t)$ is a strictly increasing (respectively, decreasing) function of $t$ on $C_+$ (respectively, on $C_-$).
\end{obs}

The proof of Theorem \ref{main} rests on two pillars: one is the characterization of elements of  best approximation \eqref{char} and the other one is the implicit function theorem (cf. \cite[Chapter III, Section 9]{E}). Markov used the orthogonality relation that yields \eqref{char} when $p=2$ (cf. \cite[Equation $2$]{M86}) together with the chain rule (cf. \cite[Equation $5$]{M86}, assuming that the zeros are implicitly defined as differentiable functions of the parameter. Kro\'o and Peherstorfer have also followed this approach in \cite{KP87}, using, in addition, the implicit function theorem to prove that the zeros are differentiable functions of the parameter. In some steps of our proof, the reader will be addressed to the corresponding step in Markov's work.  

\section{Proof of Theorem \ref{main}}\label{proof}
{\em Differentiability of the zeros}:
Let $P_{n+1}(x):=(x-x_0)\cdots(x-x_n)$, $x_j \in \re$ $(j=0,\dots, n)$. (Note that the $x_j$'s do not depend on $t$.) Define the map $\mathrm{f}:=(f_0,\dots,f_n): \mathbb{R}^{n+1} \times U \rightarrow \mathbb{R}^{n+1}$, where we have set $\mathrm{x}:=(x_0,\dots,x_n)$ and
\begin{align}\label{eq1}
f_k(\mathrm{x},t):=\int \frac{|P_{n+1}(x)|^{p} }{x-x_k}\mathrm{d}\mu(x,t).%+\sum_{i=0}^j \jmath_i(t) \frac{|P_{n+1}(y_i(t))|^{p}}{y_i(t)-x_k}.
\end{align}
For $j\not=k$ one has
\begin{align}\label{daj}
\frac{\partial f_k}{\partial x_j}(\mathrm{x},t)=&
p \int  \frac{1}{x-x_k}\frac{\partial P_{n+1}}{\partial x_j}(x) |P_{n+1}(x)|^{p-1} \text{sgn}(P_{n+1}(x)) \mathrm{d}\mu(x,t);%\\
%\nonumber  & +\sum_{i=0}^j \jmath_i(t) \frac{|P_{n+1}(y_i(t))|^{p-1}}{y_i(t)-x_k} \text{sgn}(P_{n+1}(y_i(t))) \frac{\partial P_{n+1}}{\partial x_j}(y_i(t))
\end{align}
otherwise \footnote{Cf. the denominator on the right-hand side of \cite[Equation $5$]{M86}.}
\begin{align}\label{dak}
\nonumber \frac{\partial f_k}{\partial x_k}(\mathrm{x},t)&=\int \left|\frac{P_{n+1}(x)}{x-x_k}\right|^p \frac{\partial}{\partial x_k} \left(\frac{|x-x_k|^p}{x-x_k}\right)\mathrm{d}\mu(x,t)\\
&= (1-p)\int \frac{|P_{n+1}(x)|^{p}}{(x-x_k)^{2}} \mathrm{d}\mu(x,t).%+\sum_{i=0}^j \jmath_i(t) \frac{|P_{n+1}(y_i(t))|^{p}}{(y_i(t)-x_k)^{2}}.
\end{align}
Set $\mathrm{x}(t):=(x_0(t),\dots,x_n(t))$. Fix $t_0 \in U$. From \eqref{eq1}, \eqref{daj} and \eqref{dak}, and using \eqref{char} we obtain
 $$
\mathrm{f}(\mathrm{x}(t_0),t_0)=0, \ \ \ \frac{\partial \mathrm{f}}{\partial \mathrm{x}}(\mathrm{x}(t_0),t_0)=\det \begin{pmatrix}
\dps \frac{\partial f_0}{\partial x_0}(\mathrm{x}(t_0),t_0) & &\\
& \ddots & &\\
& & \dps \frac{\partial f_n}{\partial x_n}(\mathrm{x}(t_0),t_0)&
\end{pmatrix}
\not=0.
$$ % (cf. \cite[Theorem $3.4$]{E73})
According to the implicit function theorem, under these conditions the equation $\mathrm{f}(\mathrm{s},t)=0$ has a solution $\mathrm{s}=\mathrm{x}(t)$ in a neighborhood of $(\mathrm{x}(t_0),t_0)$ that depends differentiable on $t$.

{\em Expression for the derivative of the zeros}: In view of the above result \footnote{Cf. the left-hand side of \cite[Equation $5$]{M86}.},
\begin{align*}
\frac{\mathrm{d} x_k}{\mathrm{d} t}(t)=-\frac{ \frac{\dps \partial f_k}{\dps \partial t}(\mathrm{x}(t),t)}{ \frac{\dps \partial f_k}{\dps \partial x_k}(\mathrm{x}(t),t)}.
\end{align*}
We see at once that 
\begin{align}\label{dt}
\frac{\partial f_k}{\partial t}(\mathrm{x}(t),t)=&\int \frac{|P_{n+1,p}(x,t)|^{p}}{x-x_k(t)} \frac{\partial \omega}{\partial t}(x,t)\mathrm{d}\nu(x)\\
\nonumber &+\big(\jmath'(t)+\jmath(t) y'(t) R_k(t)\big) \frac{|P_{n+1,p}(y(t),t)|^{p}}{y(t)-x_k(t)}.
\end{align}
Clearly \footnote{Cf. \cite[p. $179$]{M86}.}  
$$
\frac{1}{\omega(x_k(t),t)}\frac{\partial \omega}{\partial t}(x_k(t),t)\int \frac{|P_{n+1,p}(x,t)|^{p} }{x-x_k(t)}\mathrm{d}\mu(x,t)=0.
$$
Subtracting this from the left-hand side of \eqref{dt} yields \footnote{Cf. the numerator on the right-hand side of \cite[Equation $5$]{M86}.}
\begin{align}\label{cond}
&\frac{\partial f_k}{\partial t}(\mathrm{x}(t),t)\\
\nonumber=&\int \frac{|P_{n+1,p}(x,t)|^{p}}{x-x_k(t)}\\ 
&\left(\frac{1}{\omega(x,t)}\frac{\partial \omega}{\partial t}(x,t)-\frac{1}{\omega(x_k(t),t)}\frac{\partial \omega}{\partial t}(x_k(t),t)\right)\omega(x,t)\mathrm{d}\nu(x)\\
\nonumber &+\left(\jmath'(t)+\jmath(t) y'(t) R_k(t)-\frac{\jmath(t)}{\omega(x_k(t),t)}\frac{\partial \omega}{\partial t}(x_k(t),t)\right) \frac{|P_{n+1,p}(y(t),t)|^{p}}{y(t)-x_k(t)}.
\end{align}
It only remains to note that \footnote{Cf. \cite[p. $179$]{M86}.}
$$
\frac{1}{x-x_k(t)}\left(\frac{1}{\omega(x,t)}\frac{\partial \omega}{\partial t}(x,t)-\frac{1}{\omega(x_k(t),t)}\frac{\partial \omega}{\partial t}(x_k(t),t)\right)\geq 0.
$$
Thus
$$
\text{sgn} \left( \frac{\mathrm{d} x_k}{\mathrm{d} t}(t) \right)=\text{sgn} \left( \frac{\partial f_k}{\partial t}(\mathrm{x}(t),t) \right), 
$$
and the desired result follows from \eqref{cond}.

\section{acknowledgements}
 This work was partially supported by the Centre for Mathematics of the University of Coimbra -- UID/MAT/00324/2019, funded by the Portuguese Government through FCT/MEC and co-funded by the European Regional  Development Fund through the Partnership Agreement PT2020. FRR  is supported by the Funda\c{c}\~ao de Amparo \`a Pesquisa do Estado de Minas Gerais (FAPEMIG) under the grant PPM-00478-15.

\bibliographystyle{spbasic}      % basic style, author-year citations
\bibliography{bib}   % name your BibTeX data base

\begin{thebibliography}{10}

\bibitem{B68}
N.~Bourbaki.
\newblock {\em Elements of {M}athematics: {T}heory of {S}ets}.
\newblock Translated from the French Hermann, Publishers in Arts and Science,
  Paris; Addison-Wesley Publishing Co., Reading, Mass.-London-Don Mills, Ont.
  1968.

\bibitem{C17}
K.~Castillo.
\newblock On monotonicity of zeros of paraorthogonal polynomials on the unit
  circle.
\newblock Technical Report 17-25, Centre for Mathematics, University of
  Coimbra, 2017.

\bibitem{CR15}
K.~Castillo and F.~R. Rafaeli.
\newblock On the discrete extension of {M}arkov's theorem on monotonicity of
  zeros.
\newblock {\em Electron. Trans. Numer. Anal.}, 44:271--280, 2015.

\bibitem{C36}
J.~A. Clarkson.
\newblock Uniformly convex spaces.
\newblock {\em Trans. Amer. Math. Soc.}, 40:415--420, 1936.

\bibitem{D75}
P.~J. Davis.
\newblock {\em Interpolation and approximation. {R}epublication, with minor
  corrections, of the 1963 original, with a new preface and bibliography}.
\newblock Dover Publications, Inc., New York, 1975.

\bibitem{G71}
G.~Freud.
\newblock {\em Orthogonal polynomials}.
\newblock Pergamon Press, Oxford-New York, 1971.

\bibitem{I89}
M.~E.~H. Ismail.
\newblock Monotonicity of zeros of orthogonal polynomials.
\newblock In {\em q-{S}eries and {P}artitions (Minneapolis, MN, 1988)},
  volume~18 of {\em IMA Vol. Math. Appl.}, pages 177--190. Springer, New York,
  1988.

\bibitem{I05}
M.~E.~H. Ismail.
\newblock {\em Classical and quantum orthogonal polynomials in one variable},
  volume~98 of {\em Encyclopedia of Mathematics and Its Applications}.
\newblock Cambridge University Press, Cambridge, 2005.

\bibitem{KP87}
A.~Kro\'o and F.~Peherstorfer.
\newblock On the zeros of polynomials of minimal ${L}_p$-norm.
\newblock {\em Proc. Amer. Math. Soc.}, 101:652--656, 1987.

\bibitem{L02}
P.~Lax.
\newblock {\em Functional {A}nalysis}.
\newblock Pure and Applied Mathematics. Wiley-Interscience [John Wiley \&
  Sons], New York, 2002.

\bibitem{M86}
A.~Markov.
\newblock Sur les racines de certaines \'equations (second note).
\newblock {\em Math. Ann.}, 27:177--182, 1886.

\bibitem{S70}
I.~Singer.
\newblock {\em Best approximation in normed linear spaces by elements of linear
  subspaces}.
\newblock Translated from the Romanian by Radu Georgescu. Die Grundlehren der
  mathematischen Wissenschaften, Band 171 Publishing House of the Academy of
  the Socialist Republic of Romania, Bucharest; Springer-Verlag, New
  York-Berlin, 1970.

\bibitem{S87}
T.~J. Stieltjes.
\newblock Sur les racines de l'equation ${X}_n=0$.
\newblock {\em Acta Math.}, 9:385--400, 1887.

\bibitem{S75}
G.~Szeg\H{o}.
\newblock {\em Orthogonal polynomials}, volume~23.
\newblock Amer. Math. Soc. Coll. Publ., Amer. Math. Soc., Providence, R. I.,
  4th edition, 1975 edition, 1939.

\end{thebibliography}

\end{document}